\documentclass[preprint,3p]{elsarticle}

\usepackage{amssymb}

\usepackage{amsmath}
\usepackage{pgfplots}
\usepackage{CJKutf8}
\usepackage{setspace}
\usepackage{titlesec}
\allowdisplaybreaks[4]
\journal{arXiv}
\date{}
\begin{document}

\begin{frontmatter}



\title{Weighted Weak-type Inequalities For Fractionally Sparsely Dominated Operators}


\author[inst1]{Yanhan Chen}
\ead{chen.yanhan.67s@st.kyoto-u.ac.jp}
\affiliation[inst1]{organization={Department of Mathematics, Graduate School of Science, Kyoto University},
            postcode={Kyoto 606-8502}, 
            country={Japan}}

\begin{abstract}
In this paper, we establish quantitative weak type estimates for operators that are dominated by (fractional) sparse operators. Specifically, we derive bounds for both the restricted weak type $L^{p,1}\rightarrow L^{q,\infty}$ and the multiplier weak type, the latter of which has been previously considered by Cruz-Uribe and Sweeting \cite{CS23}. These estimates provide a precise quantification of the mapping properties of the considered operators, extending and refining the existing theory.
\end{abstract}




\end{frontmatter}


\section{Introduction}
\label{sec:sample1}
Over the past decade, sparse operators have garnered significant attention within the framework of weighted theory, particularly following the work of Moen \cite{M12}, who provided an alternative proof of the renowned $A_{2}$ conjecture utilizing sparse methods after the initial proof by Hytönen \cite{H12}. Therein, Moen \cite{M12} has derived a sharp weighted estimate for the operator $T$, which satisfies the strong pointwise sparse domination 
$$|Tf(x)|\lesssim\sum_{Q\in\mathcal{S}}\langle f\rangle_{Q}\chi_{Q}(x):=\mathcal{A}_{\mathcal{S}}f(x),$$
for some sparse family $\mathcal{S}$, where $\langle f\rangle_{Q}:={|Q|}^{-1}\int_{Q}|f|dx$. To be precisely, for $A_{p}$ weight $w$ and $1<p<\infty$, it holds that
$$\|\mathcal{A}_{\mathcal{S}}f\|_{L^{p}(w)}\leq c_{p}[w]_{A_{p}}^{\text{max}\left(1,\frac{p^{'}}{p}\right)}\|f\|_{L^{p}(w)}$$
with a positive constant $c_{p}$.
Later it was shown by Bernicot, Frey and Petermichl \cite{BFP16} that the method of sparse operators possesses a significantly broader scope than Calder\'{o}n-Zygmund operators. For instance, one can consider the Bochner-Riesz operator $B_{\lambda}f=\mathcal{F}^{-1}[(1-|\cdot|^{2})^{\lambda}_{+}\mathcal{F}f]$ for $\lambda>0$, or spherical maximal operator $M_{S}f:=\sup_{t>0}\left|\int_{S^{n-1}}f(\cdot-ty)d\sigma(y)\right|$. Generally, such operators can not be pointwisely bounded by sparse operators as before. Nevertheless, they do possess a sparse domination characteristic that actually derives the quantitative weighted bounds, see for instance \cite{KL18,LMR19,L19}. Such operators satisfy the $L^{p}$ boundedness for only restricted $p\in(p_{0},q_{0})$, Bernicot, Frey and Petermichl \cite{BFP16} gave a precise description of the common characteristic of those operators with a weak sparse domination as following
$$|\langle T f, g\rangle| \lesssim\sum_{Q \in\mathcal{S}}\langle f\rangle_{p_0, Q}\langle g\rangle_{q_0^{\prime}, Q}|Q|$$
for any $g\in L^{\infty}$ with compact support. This definition covers the strong pointwise domination within the case $p_{0}=1$ and $q_{0}=\infty$. Furthermore, they demonstrated that such bilinear sparse domination leads to optimal quantitative $A_{p}-RH_{s}$ estimates for $p_{0}<p<\infty$.\\
\indent While for Calder\'{o}n-Zygmund singular integral $T$, improvements can be made to the weak-type bounds inherited from the $A_{2}$ strong-type bound. For example, according to the work of Hyt\"{o}nen, Lacey, Martikainen, Orponen, Reguera, Sawyer, Uriarte-Tuero \cite{HLM12}, one has 
$$\|Tf\|_{L^{1,\infty}(w)}\lesssim(1+\text{log}[w]_{A_{1}})[w]_{A_{1}}\|f\|_{L^{1}(w)}$$
for $w\in A_{1}$. Hyt\"{o}nen and P\'{e}rez \cite{HP13} found that the bounds could be further improved when introducing the smaller Fujii-Wilson constant $[w]_{A_{\infty}}$, there the weak type constant could be replaced by $(1+\text{log}[w]_{A_{\infty}})[w]_{A_{1}}$. As for general (bilinear) sparsely dominated operator, Frey and Nieraeth \cite{FN19} proved weighted weak type $(p_{0},p_{0})$ boundedness $L^{p_{0}}\rightarrow L^{p_{0},\infty}$ with quantitative mixed $A_{1}-A_{\infty}$ estimates.\\
\indent Significant attention has also been devoted to the study of some fractional-type operators. Cruz-Uribe and Moen  \cite{CM13} varified that Riesz potentials $I_{\alpha}$ (or singular integral) could be pointwisely bounded by some fractional sparse operator that
$$|I_{\alpha}f(x)|\lesssim\sum_{i=1}^{N}\sum_{Q\in\mathcal{S}_{i}}|Q|^{\frac{\alpha}{n}}\langle f\rangle_{Q}\chi_{Q}(x).$$
Fractional maximal operator $M_{\alpha}f:=\sup_{Q}\langle f\rangle_{\alpha,1,Q}\chi_{Q}$ is another example of this pointwise sparsely domination. Moen \cite{M12} has also considered the weighted $A_{2}$-inequality of such fractional sparse operator $\mathcal{A}_{\mathcal{S}}^{\alpha}f(x):=\sum_{Q\in\mathcal{S}}\langle f\rangle_{\alpha,1,Q}\chi_{Q}(x)$, where $\langle f\rangle_{\alpha,p,Q}:=\left({|Q|^{-1+{\alpha}/{n}}}\int_{Q}|f|^{p}dx\right)^{{1}/{p}}$. They showed that for $0<\alpha<n$, $1<p<\alpha/n$, $1/q=1/p-\alpha/n$ satisfy $\text{min}(p^{'}/q,q/p^{'})\leq1-\alpha/n$ and $w\in A_{p,q}$, the following estimate holds
$$\|\mathcal{A}_{\mathcal{S}}^{\alpha}f\|_{L^{q}(w^{q})}\leq [w]_{A_{p,q}}^{\left(1-\frac{\alpha}{n}\right)\text{max}\left(1,\frac{p^{'}}{q}\right)}c_{p,\alpha}\|f\|_{L^{p}(w^{p})}$$
with a positive constant $c_{p,\alpha}$. While there is somewhat unnatural that $\text{min}(p^{'}/q,q/p^{'})\leq1-\alpha/n$ in the assumption, thus they do not obtain the full range of $p,q$ for even Riesz potentials $I_{\alpha}$ or fractional maximal operator $M_{\alpha}$. Recently, when Lerner, Lorist and Ombrosi \cite{LLO24} considered Bloom weighted estimate for sparse forms associated to commutators, they have introduced such general fractional bilinear sparse domination for operator $T$
\begin{equation}
|\langle T f, g\rangle| \leq C\sum_{j=1}^{N}\sum_{Q \in\mathcal{S}_{j}}\langle f\rangle_{p_0, Q}\langle g\rangle_{\alpha,q_0^{\prime}, Q}|Q|.\tag{1.1}
\end{equation}
Such definition covers the notation of pointwise domination by fractional sparse operator $\mathcal{A}_{\mathcal{S}}$ within the case $p_{0}=1$ and $q_{0}=\infty$. They also proved a more general $L^{p}(w^{p})\rightarrow L^{q}(w^{q})$ weighted inequality with an entire new approach that technically removing the unnatural assumption.\\
\indent The goal of the current work is to establish quantitative weighted weak estimate including the endpoint $p=p_{0}$ for operators satisfy (1.1). We first introduce a bit weaker $A_{p}$ weights ($A^{\mathcal{R}_{p}}$) into a tough restricted weak type argument for fractional sparse operator, where a Calderón-Zygmund decomposition is applicable. Then we derive bounds for multiplier weak type bounds follows with the idea by Cruz-Uribe and Sweeting \cite{CS23}.
\subsection{Main results}
\indent Let $0\leq\alpha<n$ and $1\leq p_{0}<q_{0}\leq\infty$ satisfy ${1}/{p_{0}}-{1}/{q_{0}}>{\alpha}/{(nq_{0}^{'})}$. We consider an operator $T\in S^{\alpha}(p_{0},q_{0})$, which means the bound $(1.1)$ holds for $T$, the precise definition could be found in Section 2. When $\alpha=0$, it has been shown by Bernicot, Frey and Petermichl \cite{BFP16} that for $T\in S^{0}(p_{0},q_{0})$, $p_{0}<p<q_{0}$ and $w\in A_{p/p_{0}}\cap RH_{(q_{0}/p)^{'}}$,
\begin{equation}
    ||T||_{L^{p}(w)\rightarrow L^{p}(w)}\lesssim\left([w]_{A_{\frac{p}{p_{0}}}}[w]_{RH_{\left(\frac{q_{0}}{p}\right)^{'}}}\right)^{\text{max}\left(\frac{1}{p-p_{0}},\frac{q_{0}-1}{q_{0}-p}\right)}.\tag{1.2}
\end{equation}
Later, the two-weight strong type estimates $T(\cdot\sigma): L^{p}(\sigma)\rightarrow L^{q}(w)$ for the particular case when $q_{0}=\infty$ are encompassed within the research conducted by Fackler and Hytönen \cite{FH17}. Recently, Lerner, Lorist and Ombrosi \cite{LLO24} have proved the general bounds for $0\leq\alpha<n$ with the test condition method introduced by Li \cite{L17}. With a similar method, their result could be slightly strengthened as follows.\\
~\\
\noindent\textbf{Theorem A} (Lerner, Lorist, Ombrosi \cite{LLO24}, modified): \textit{Let $1\leq p_{0}<p\leq q<q_{0}\leq\infty$, ${1}/{p}-{1}/{q}={\alpha}/{nq_{0}^{'}}\ (0\leq\alpha<n)$, $T\in S^{\alpha}(p_{0},q_{0})$ and $w^{q}\in A_{\left(\frac{1}{p_{0}}-\frac{1}{p}\right)q+1}\cap RH_{\left(\frac{q_{0}}{q}\right)^{'}}$. Then
    \begin{equation}
        ||T||_{L^{p}(w^{p})\rightarrow L^{q}(w^{q})}\lesssim\left([w^{q}]_{A_{\left(\frac{1}{p_{0}}-\frac{1}{p}\right)q+1}}[w^{q}]_{RH_{\left(\frac{q_{0}}{q}\right)^{'}}}\right)^{\theta}\tag{1.3}
    \end{equation}
with} $$\theta=\text{max}\left\{\left(\frac{q_{0}}{q}\right)^{'}\frac{1-\frac{\alpha}{n}}{q_{0}^{'}},\frac{\frac{1}{p_{0}}-\frac{\alpha}{nq_{0}^{'}}}{q\left(\frac{1}{p_{0}}-\frac{1}{p}\right)}\right\}.$$\\
\indent In the work of Frey and Nieraeth \cite{FN19}, it was demonstrated that under the condition $\alpha=0$, the operator $T$ as defined in (1.1) adheres to a weak type inequality precisely at the endpoint $p=p_{0}$, taking the form $\|T\|_{L^{p}(w) \rightarrow L^{p, \infty}(w)}<\infty$. This inequality is intricately linked to the $A_{1}-A_{\infty}$ constant associated with the weight function $w$.\\
\indent Recently, Fay, Rey and Škreb \cite{FRS24} considered the restricted weak type estimate for sparse operator $\mathcal{A}_{\mathcal{S}}f:=\sum_{Q\in\mathcal{S}}\langle f\rangle_{Q}\chi_{Q}$. 
They found the exact Bellman function associated to the level-sets of sparse operators acting on characteristic functions. In the general cases including the endpoint, we establish the following restricted weak type estimate.\\

\noindent \textbf{Theorem 1.1}: \textit{Let $1\leq p_{0}\leq p\leq q<q_{0}\leq\infty$, $1/p-1/q=\alpha/(nq_{0}^{'})\ (0\leq\alpha<n)$, $0\leq\eta\leq1$, $T\in S^{\alpha}(p_{0},q_{0})$ and $w^{p_{0}}\in A^{\mathcal{R}}_{p/p_{0},q/p_{0}}$, $w^{q}\in RH_{\left(q_{0}/q\right)^{'}}$. Then}
$$\|T\|_{L^{p, 1}(w^{p}) \rightarrow L^{q, \infty}(w^{q})}\lesssim[w^{q}]_{A_{\frac{p}{p_{0}},\frac{q}{p_{0}}}^{\mathcal{R}}}^{\frac{1}{p_{0}}}[w^{q}]_{RH_{\left(\frac{q_{0}}{q}\right)^{'}}}^{\left(\frac{q_{0}}{q}\right)^{'}+\frac{1}{q}}D^{\eta}_{w^{q}}\left(D_{w^{q}}^{2^{n}}\right)^{\frac{1}{q}},$$
\textit{where $[w^{q}]_{RH_{1}}=1$,}
$$D_{w^{q}}^{\eta}=\sup\left\{\frac{w^{q}(Q)}{w^{q}(E)}:\ Q\text{ is a shifted dyadic cube},\ E\subseteq Q,\ |E|\geq\eta|Q|\right\}\lesssim[w^{q}]_{A_{\frac{p}{p_{0}},\frac{q}{p_{0}}}^{\mathcal{R}}}^{\frac{q}{p_{0}}}.$$\\
\indent According to the work of Duoandikoetxea and Mart\'{i}n-Reyes \cite{DM16}, there exists a representation of $D^{\eta}_{w^{q}}$ that is solely dependent on the reverse Hölder constant $[w^{q}]_{RH_{(q_{0}/q)^{'}}}$, and in the special case when $q_{0}=\infty$, this constant is replaced by the Fujii-Wilson constant $[w^{q}]_{A_{\infty}}$. However, it is important to note that the involved estimate is rather crude, indicating that the representation does not necessarily yield an optimal constant.\\ 
\indent In the case $p_{0}=1$, Kokilashivili \cite{K87} presented an exemplary fractional maximal operator defined as
$M_{\alpha}f(x):=\sup_{Q}\langle f\rangle_{\alpha,1,Q}\chi_{Q}(x)$, which belongs to the class $S^{\alpha}(1,\infty)$. This example (together with Proposition 2.5) shows the necessity of the weight class $A^{\mathcal{R}}_{p,q}$ and the sharpness of the associated constant $[w^{q}]_{A^{\mathcal{R}}_{p,q}}$. Furthermore, such sharpness in both weight class and constant could be extended to general $p_{0}\geq1$ if we consider the $p_{0}$-version of $M_{\alpha}$. While it remains uncertain whether the sharpness of reverse Hölder class $RH_{(q_{0}/q)^{'}}$ and its associated constant hold true for the case when $q_{0}<\infty$.\\
\indent Besides, we conjecture that the doubling constant $D^{\eta}_{w^{q}}$ is not a prerequisite in this context, and the prospect of investigating the restricted weak type boundedness of such operators without relying on the those constant is also worth considering.\\
~\\
\textbf{Remark 1.2:} The Marcinkiewicz’s interpolation theorem could be extended to an off-diagonal version: a pair of restricted weak type estimates $L^{p_{i},1}(u)\rightarrow L^{q_{i},\infty}(v)$ ($i=1,2$) are powerful enough to often imply quantitative estimates on intermediate Lorentz spaces $L^{p,r}(u)\rightarrow L^{q,r}(v)$ for any $1\leq r\leq\infty$ \cite{SM93}.\\
~\\
\indent We also consider multiplier weak type bound. Typically, the weighted strong type inequality
\begin{equation}   
||Tf||_{L^{q}(w^{q})}\lesssim||f||_{L^{p}(w^{p})}\tag{1.4}
\end{equation}
implies the weak type estimate 
$$||Tf||_{L^{q,\infty}(w^{q})}\lesssim||f||_{L^{p}(w^{p})}$$
according to Chebyshev's inequality. While if we treat $w$ as a multiplier instead of weight, the inequality $(1.4)$ leads to the multiplier weak type bound
$$||wT(w^{-1}f)||_{L^{q,\infty}}\lesssim||f||_{L^{p}}$$
Cruz-Uribe and Sweeting \cite{CS23} proved such inequality for the maximal operator and Calder\'{o}n-Zygmund singular integrals with the sparse method. As an extension of their result, we provide that\\
~\\
\textbf{Theorem 1.3}: \textit{Let $1\leq p_{0}\leq p\leq q<q_{0}\leq\infty$, $1/p-1/q=\alpha/(nq_{0}^{'})\ (0\leq\alpha<n)$, $T\in S^{\alpha}(p_{0},q_{0})$ and $w^{q}\in A_{\left(1
/p_{0}-1/p\right)q+1}\cap RH_{\left(q_{0}/q\right)^{'}}$. Then}

$$||wT(w^{-1}\cdot)||_{L^{p}\rightarrow L^{q,\infty}}\leq
\begin{cases}
    [w^{q}]_{A_{\left(\frac{1}{p_{0}}-\frac{1}{p}\right)q+1}}^{\frac{1}{q}}[w^{q}]_{A_{\infty}}\quad\quad\quad q_{0}=\infty\\
    [w^{q}]_{A_{\left(\frac{1}{p_{0}}-\frac{1}{p}\right)q+1}}^{\frac{1}{q}}[w^{q}]_{RH_{\left(\frac{q_{0}}{q}\right)^{'}}}^{\left(\frac{q_{0}}{q}\right)^{'}+\frac{2}{q}}\quad 1<q_{0}<\infty
\end{cases}.$$
\indent As we mentioned previously, multiplier weak type inequalities can be derived from strong type inequalities. Theorem A leads to the quantitative estimate
$$||wT(w^{-1}\cdot)||_{L^{p}\rightarrow L^{q,\infty}}\lesssim\left([w^{q}]_{A_{\left(\frac{1}{p_{0}}-\frac{1}{p}\right)q+1}}[w^{q}]_{RH_{\left(\frac{q_{0}}{q}\right)^{'}}}\right)^{\theta},$$
where $\theta$ is the same as defined above. It could be simply check that $1/q<\theta$, which means for non-endpoint case $p>p_{0}$ our estimate is still new. For the ordinary case $T=\mathcal{A}_{\mathcal{S}}^{\alpha}$, Cruz-Uribe and Sweeting \cite{CS23} gave the same constant $[w^{q}]_{A_{{q}/{p^{'}}+1}}^{{1}/{q}}[w^{q}]_{A_{\infty}}$. While in the case $q_{0}<\infty$, our constant is much more larger, we believe it could be greatly deduced.\\

\noindent\textbf{Remark 1.4:} It still remains an open question, even for maximal operator or singular integral, to ascertain the necessary and sufficient conditions on w for the multiplier weak-type inequality to hold true.
\section{Preliminaries}
\subsection{The setting}
\indent In this paper we only consider the Euclidean space $\mathbf{R}^{n}$ equipped with Lebesgue measure. While the theorems and discussion below could be transplanted to general Borel measure $\mu$ that is finite on compact sets and strictly positive on non-empty open set with doubling condition, i.e. there is a constant $C>0$ such that
$$\mu(2B)\leq C\mu(B)$$
for any balls $B$ and the ball with the same center as $B$ whose radius is twice of the radius of $B$ is denoted by $2B$. We will use $a\lesssim b$ to say that there exists a constant $C$, which is independent of the important parameters, such that $a\leq Cb$. Moreover, we write $a\sim b$ if $a\lesssim b$ and $b\lesssim a$.\\
\indent For any measurable set $E$, let us denote the Lebesgue measure of $E$ by $|E|$, and for any non-negative weight $w$, denote $w(E):=\int_{E}wdx$, and write $\chi_{E}$ as the characteristic function of the set $E$. For any $1\leq p<\infty$ and $0\leq\alpha<n$, we will write
$$\langle f\rangle_{\alpha,p,E}:=\left(\frac{1}{|E|^{1-\frac{\alpha}{n}}}\int_{E}|f|^{p}dx\right)^{\frac{1}{p}},$$
when $\alpha=0$ we simply write $\langle f\rangle_{p,E}=\langle f\rangle_{0,p,E}$ and $\langle f\rangle_{E}=\langle f\rangle_{1,E}$. We write $\langle f,g\rangle:=\int fgdx$, and define $p^{'}:={p}/{p-1}$ for any $1\leq p\leq\infty$.\\
\indent A collection $\mathcal{D}$ of cubes in $\mathbf{R}^{n}$ is called a \textit{dyadic grid} if for each cube in $\mathcal{D}$, the sides of the cube are parallel to axis and $\mathcal{D}$ satisfies the following properties:
\begin{align}
&\text{(1) For any $Q\in\mathcal{D}$, its side length $\ell(Q)$ is of the form $2^{k}$, for some $k\in\mathbf{Z}$.}\notag\\
&\text{(2) For all $Q$ and $R\in\mathcal{D}$, $Q\cap R\in\left\{\varnothing,Q,R\right\}$.}\notag\\
&\text{(3) The cubes of a fixed side length $2^{k}$ form a partition of $\mathbf{R}^{n}$. }\notag
\end{align}
The standard dyadic grid in $\mathbf{R}^{n}$ consists of cubes $2^{-k}([0,1)^{n}+j),$ with $k\in\mathbf{Z}$ and $j\in\mathbf{Z}^{n}$. The shifted dyadic cubes are defined by 
$$\mathcal{D}^{a}:=\left\{2^{-k}\left([0,1)^{n}+m+(-1)^{k}\frac{a}{3}\right);\ k\in\mathbb{Z},m\in\mathbb{Z}^{n}\right\},$$
with $a\in\left\{0,1,2\right\}^{n}$.\\
\indent For a dyadic grid $\mathcal{D}$, a collection of cubes $\mathcal{S}\subseteq\mathcal{D}$ is said to be \textit{$\eta$-sparse} for some $0<\eta\leq 1$, if there is a pairwise disjoint collection $(E_Q)_{Q\in S}$, so that $E_Q\subseteq Q$, $|E_Q|\geq\eta|Q|$.\\
~\\
\textbf{Definition 2.1:}\textit{ Let $T$ be a (sub)linear operator, initially defined on $C_{c}^{\infty}$, with the following property: There are $1\leq p_0<q_0\leq\infty$, a positive integer $N$, $0<\eta\leq 1$ and $0<\alpha<n$, so that for each $f,g\in C_{c}^{\infty}$ (which means $f,g$ are differentiable at any order and have compact support), there exist $\eta$-sparse collections $\mathcal{S}_{j}(j=1,2,...,N)$, so that
\begin{equation}
|\langle T f, g\rangle| \lesssim\sum_{j=1}^{N}\sum_{Q \in\mathcal{S}_{j}}\langle f\rangle_{p_0, Q}\langle g\rangle_{\alpha,q_0^{\prime}, Q}|Q|.\notag
\end{equation}
Then, we will say $T$ is a fractionally sparsely dominated operator, written as $T\in S^{\alpha}(p_0,q_0)$. When $\alpha=0$, we simply write $S^{0}(p_{0},q_{0})=S(p_{0},q_{0}).$}\\
~\\
\textbf{Remark 2.2: }Such operator coincides with the definition of the sparsely dominated operator in \cite{BFP16} when $\alpha=0$ (non-fraction case), which includes the Hardy-Littlewood maximal operator and the Calder\'{o}n-Zygmund singular integral. In the case $q_{0}=\infty$ and $p_{0}=1$, if we assume $TF$ to be integrable, this becomes the fractional sparse operator $A_{\mathcal{S}}^{\alpha}$ defined as
$$\mathcal{A}_{\mathcal{S}}^{\alpha}f:=\sum_{Q\in\mathcal{S}}\left(\frac{1}{|Q|^{1-\alpha/n}}\int_{Q}|f|dx\right)\chi_{Q}.$$
Those classes of operators were proven to enjoy good mapping properties in weighted $L^{p}$ spaces or relative Lorentz spaces for $p_0\leq p < q_0$.\\
~\\
\textbf{Remark 2.3: }There is a wealth of examples of fractionally sparsely dominated operators, for instance \cite{BFP16,CCPO17,L19} and references therein. The applications in PDE are also worthy of attention. Recently, Saari, Wang and Wei \cite{SWW24} have shown the local solution to the divergence form elliptic equation: $\text{div}\ a(x, \nabla u(x))=\text{div}\ F(x)+f(x)$, where $a$ is an elliptic coefficient, could be divided into two parts with each one satisfies a local bound as in (1.1) in some sense. Our general theory for (fractional) sparse operator leads to some quantitative weighted estimates of $u$.
\subsection{Lorentz space}
\indent Let us recall the definition of Lebesgue and Lorentz spaces. For $1\leq p<\infty$ and an arbitrary measure space $(X,u)$, let us denote $L^{p,1}(u)$ the Lorentz space of $u$-measurable functions 
$$||f||_{L^{p,1}}(u):=p\int_{0}^{\infty}\lambda^{u}_{f}(y)^{\frac{1}{p}}dy=\int_{0}^{\infty}f_{u}^{*}(t)t^{\frac{1}{p}}\frac{dt}{t}<\infty,$$
where $\lambda^{u}_{f}$ denotes the distribution function of $f$ and $f^{*}_{u}$ is the decreasing rearrangement of $f$ with respect to $u$
$$\lambda^{u}_{f}:=u(\{x\in X:|f(x)|>t\})\ \text{ and }\ f_{u}^{*}(t):=\inf\{y>0:\lambda^{u}_{f}(y)\leq t\}.$$
\noindent We denote the Lebesgue space of $u-$measurable functions by $L^{p}(u)$ such that
$$||f||_{L^{p}(u)}:=\left(\int_{X}|f|^{p}du\right)^{\frac{1}{p}}<\infty,$$
and $L^{p,\infty}(u)$ is the Lorentz space of $u$-measurable functions such that
$$||f||_{L^{p,\infty}(u)}:=\sup_{y>0}y\lambda^{u}_{f}(y)^{\frac{1}{p}}=\sup_{t>0}t^{\frac{1}{p}}f^{*}_{u}(t)<\infty.$$
It is well-known that $L^{p,1}(u)\hookrightarrow L^{p}(u)\hookrightarrow L^{p,\infty}(u)$ for $1\leq p<\infty$.
\subsection{Weight class}
\indent We identify a weight $w$ with a Borel measure by setting $w(E):=\int_{E}wdx$ for all measurable sets $E\subseteq\mathbf{R}^{n}$. We first define the Muckenhoupt weight \cite{M72}. For $1\leq p<\infty$ we say that $w\in A_{p}$ if
$$[w]_{A_{p}}:=\sup_{Q}\left(\frac{1}{|Q|}\int_{Q}wdx\right)\left(\frac{1}{|Q|}\int_{Q}w^{1-p^{'}}dx\right)^{p-1}<\infty,$$
where for $p=1$ we use the limiting interpretation $(\int_{Q}w^{1-p^{'}}dx/|Q|)^{p-1}=(\text{essinf}_{Q}\ w)^{-1}$. Define the class of weights
$$A_{\infty}:=\bigcup_{p\geq1}A_{p}.$$
Then a weight $w\in A_{\infty}$ if and only if
$$[w]_{A_{\infty}}:=\sup_{Q}\frac{1}{w(Q)}\int_{Q}M(w\chi_{Q})dx<\infty,$$
where $M$ denotes the Hardy-Littlewood maximal operator
$$Mf(x)=\sup_{Q}\langle f\rangle_{Q}\chi_{Q}(x).$$
This quantity is referred to as the Fujii-Wilson $A_{\infty}$ constant \cite{F77,W87}.\\
\indent For $1\leq p<\infty$, following, Hunt and Kurtz \cite{CHK82}, we say that $w\in A_{p}^{\mathcal{R}}$  if 
\begin{align*}
[w]_{A^{\mathcal{R}}_{p}}&:=\sup_{Q}\sup_{E\subseteq Q}\frac{|E|}{|Q|}\left(\frac{w(Q)}{w(E)}\right)^{\frac{1}{p}}\\
&\sim\sup_{Q}w(Q)^{\frac{1}{p}}\frac{||\chi_{Q}w^{-1}||_{L^{p^{'},\infty}(w)}}{|Q|}=:[w]^{'}_{A^{\mathcal{R}}_{p}}<\infty.   
\end{align*}
This class of weights is consistently taken into account when examining the restricted weak-type estimate $L^{p,1}(w)\rightarrow L^{p\infty}(w)$. Analogous to the first definition presented, in this paper, we further consider a $p-q$ version restricted weak type weight.\\ 
~\\
\textbf{Definition 2.4:} For $1\leq p\leq q<\infty$, $1/p-1/q=\alpha/n$ $(0\leq\alpha<n)$, we say that $w\in A_{p,q}^{\mathcal{R}}$ if
$$[w]_{A_{p,q}^{\mathcal{R}}}:=\sup_{Q}\sup_{E\subseteq Q}\frac{|E|}{|Q|^{1-\frac{\alpha}{n}}}\frac{(w^{q}(Q))^{\frac{1}{q}}}{(w^{p}(E))^{\frac{1}{p}}}<\infty.$$
This definition aligns with Koilashivili's definition, as outlined in \cite{K87} (espaecially take $w,\phi,v$ as $w^{q},1,w^{p}$) in the following sense:\\
~\\
\textbf{Proposition 2.5:} For $1\leq p\leq q<\infty$, $1/p-1/q=\alpha/n$ $(0\leq\alpha<n)$, define 
$$[w]_{A_{p,q}^{\mathcal{R}}}^{'}:=\sup_{Q}w^{q}(Q)^{\frac{1}{q}}\|\chi_{Q}w^{-p}\|_{L^{p^{'},\infty}(w^{p})}|Q|^{\frac{\alpha}{n}-1}.$$
Then $[w]_{A_{p,q}^{\mathcal{R}}}^{'}\sim[w]_{A_{p,q}^{\mathcal{R}}}$.\\
\textit{Proof: }The proof is concise. Firstly, for any $E\subseteq Q$, apply H\"older's inequality, it holds $$|E|=\int_{E}w^{-p}w^{p}dx\lesssim\|\chi_{E}w^{-p}\|_{L^{p^{'},\infty}(w^{p})}\|\chi_{E}\|_{L^{p,1}(w^{p})}\leq[w]_{A_{p,q}^{\mathcal{R}}}|Q|^{1-\frac{\alpha}{n}}\frac{w^{p}(E)^{\frac{1}{p}}}{w^{q}(Q)^{\frac{1}{q}}},$$ thus $[w]_{A_{p,q}^{\mathcal{R}}}\lesssim[w]_{A_{p,q}^{\mathcal{R}}}^{'}$. Conversely, for fixed cube $Q$ and $y>0$, let $E_{y}=\{x\in Q:\ w^{-p}(x)>y\}$, then 
$$yw^{p}(E_{y})=\int_{E_{y}}yw^{p}dx\leq\int_{E_{y}}w^{-p}w^{p}dx=|E_{y}|\leq[w]_{A_{p,q}^{\mathcal{R}}}|Q|^{1-\frac{\alpha}{n}}\frac{w^{p}(E_{y})^{\frac{1}{p}}}{w^{q}(Q)^{\frac{1}{q}}},$$
which leads to $[w]_{A_{p,q}^{\mathcal{R}}}^{'}\leq[w]_{A_{p,q}^{\mathcal{R}}}$, and hence $[w]_{A_{p,q}^{\mathcal{R}}}^{'}\sim[w]_{A_{p,q}^{\mathcal{R}}}$. $\square$\\
\indent For $1<s\leq\infty$, we follow the general notation by Coifman and Fefferman \cite{CF74}, say that $w\in RH_{s}$ if
$$[w]_{RH_{s}}:=\sup_{Q}\frac{\langle w\rangle_{s,Q}}{\langle w\rangle_{Q}}<\infty.$$
For $s=1$ we will use the interpretation $RH_{1}=A_{\infty}$ and let $[w]_{RH_{1}}:=1.$\\
\indent We provide some facts about the those classes. \\
~\\
\textbf{Proposition 2.6}: (1) \textit{Given $w\in RH_{s}\ (1<s<\infty)$, there exists $c>0$ only depend on dimension $n$ such that if $v=s+(s-1)(cs[w]^{s}_{RH_{s}})^{-1}$, then $w\in RH_{v}$ and $[w]_{RH_{v}}\lesssim[w]_{RH_{s}}$.}\\
(2) \textit{Given $w\in A_{\infty}$, there exists $d>0$ only depend on dimension $n$ such that if $v=1+d[w]_{A_{\infty}}^{-1}$, then $w\in RH_{v}$ and $[w]_{RH_{v}}\leq2.$}\\
\indent Gehring gave a proof of (1) in \cite{G73}, for (2) we refer to Hyt\"onen, P\'erez and Rela's work \cite{HPR12}.
\subsection{Maximal operators}
\indent Given a Borel measure $u$ on $\mathbf{R}^{n}$ and a dyadic grid $\mathcal{D}$, we define the dyadic fractional maximal operator
$$\mathcal{M}^{\mathcal{D}}_{\alpha,u}f(x):=\sup_{Q\in\mathcal{D}}\frac{1}{u(Q)^{1-\frac{\alpha}{n}}}\int_{Q}|f|du\chi_{Q}(x),\quad0\leq\alpha<n.$$
When $\alpha=0$ we simply write $\mathcal{M}^{\mathcal{D}}_{u}=\mathcal{M}^{\mathcal{D}}_{0,u}$. we have the following property concerning the boundedness of $\mathcal{M}^{\mathcal{D}}_{\alpha,u}.$\\
~\\
\textbf{Proposition 2.7} (Moen, \cite{M12}): \textit{If $0\leq\alpha<n$, $1<p\leq n/\alpha$ and $1/p-1/q=\alpha/n$, then
$$||\mathcal{M}^{\mathcal{D}}_{\alpha,u}||_{L^{q}(u)}\leq\left(1+\frac{p^{'}}{q}\right)^{1-\frac{\alpha}{n}}||f||_{L^{p}(u)}.$$
Further at the endpoint $p=1$, $q_{0}=n/(n-\alpha)$, we have $||\mathcal{M}^{\mathcal{D}}_{\alpha,u}||_{L^{1,\infty}(u)}\leq||f||_{L^{q_{0}}(u)}.$}
\section{Proofs of the main theorems}
The sum on the right-hand side of $(1.1)$ can be split into $N$ sums by considering different dyadic grids, thus in this section these proofs only consider a single dyadic grid $\mathcal{D}$ and the $\eta$ sparse collection $\mathcal{S}\subseteq\mathcal{D}$. Besides, with a similar method as Lacey and Mena \cite{LM16}, the existence of the universal fractional sparse collection could be verified. Hence, we can further assume $\mathcal{D}=\mathcal{D}^{a}$ for $a\in\left\{0,1,2\right\}^{n}$.
\subsection{Proof of Theorem 1.1}
\indent As been well-known (see for example Stein's literature \cite{SM93}), we shall only consider that case $f=\chi_F$, where $F$ is a measurable set in $\mathbb{R}^n$, and assume that $w^{p}(F)=1$. To prove the theorem, we will use the equivalence (Grafakos, \cite{G08})
\begin{equation}    
    ||T(f)||_{L^{q,\infty}(w^{q})}\sim\sup\limits_{0<w^{q}(G)<\infty}\inf_{\substack{G^{'}\subseteq G\\w^{q}(G^{'})\geq\frac{1}{2}w^{q}(G)}}w^{q}(G^{'})^{-1+\frac{1}{q}}|\langle Tf, \chi_{G^{'}}w^{q}\rangle|.\tag{3.1}
\end{equation}
Let $G\subseteq\mathbf{R}^{n}$ with $0<w^{q}(G)<\infty$. Let 
\begin{equation}
    \Omega:=\left\{\mathcal{M}^{\mathcal{D}}_{w^{q}}(fw^{p-q})>\frac{2}{w^{q}(G)}\right\}\supseteq\left\{fw^{p-q}>\frac{2}{w^{q}(G)}\right\},\tag{3.2}
\end{equation}
where the inclusion is up to a set of measure zero.\\
\underline{The case $\Omega\neq\varnothing$ and $q_{0}<\infty$.} Let $\mathcal{P}$ be the family of those maximal cubes $P\in\mathcal{D}$ that 
$$\frac{1}{w^{q}(P)}\int_{P}fw^{p}dx>\frac{2}{w^{q}(G)}.$$
Then $\bigcup_{P\in\mathcal{P}} P=\Omega$, according to Proposition 2.7 we have
$$w^{q}(\Omega)\leq\frac{w^{q}(G)}{2}\int fw^{p-q}w^{q}dx=\frac{w^{q}(G)}{2}.$$
Choose $G^{'}:=G\cap\Omega^{c}$, then $w^{q}(G^{'})\geq1/2w^{q}(G)$. Now, since $fw^{p-q}=\chi_{F}w^{p-q}\in L_{loc}^{1}(w^{q})$, we can consider the Calder\'{o}n-Zygmund decomposition : $fw^{p-q}=g+b$, where\\
\begin{align*}
&g=\sum_{P\in\mathcal{P}}\left(\frac{1}{w^{q}(P)}\int_{P}fw^{p}dx\right)\chi_{P}+fw^{p-q}\chi_{\Omega^{c}}\ \ \text{and}\\
&b=\sum_{P\in\mathcal{P}}\left(fw^{p-q}-\frac{1}{w^{q}(P)}\int_{P}fw^{p}dx\right)\chi_{P}.
\end{align*}
Here $||g||_{L^{1}(w^{q})}=||f||_{L^{1}(w^{p})}=1$ and $\langle bw^{q}\rangle_{P}=0$ for any $P\in\mathcal{P}$. Moreover, it holds 
$(w^{q}(P))^{-1}\int_{P}fw^{p}dx>2(w^{q}(G))^{-1}$ and $(w(\tilde{P}))^{-1}\int_{\tilde{P}}fw^{p}dx\leq2(w^{q}(G))^{-1}$, where $\tilde{P}$ denotes the minimal cube in $\mathcal{D}$ that strictly contains $P$, thus $||g||_{L^{\infty}}\lesssim D^{2^{n}}_{w^{q}}w^{q}(G)^{-1}$. By the definition of $A^{\mathcal{R}}_{p,q}$ weight,

\begin{align*}
\langle f\rangle_{p_0, Q}&=\left(\frac{|F \cap Q|}{|Q|}\right)^{\frac{1}{p_0}}=\frac{1}{|Q|^{\frac{\alpha}{nq_{0}^{'}}}}\left(\frac{|F\cap Q|}{|Q|^{1-\frac{p_{0}\alpha}{nq_{0}^{'}}}}\right)^{\frac{1}{p_{0}}}\\
&\leq[w^{q}]_{A_{\frac{p}{p_{0}},\frac{q}{p_{0}}}^{\mathcal{R}}}^{\frac{1}{p_{0}}}\frac{1}{|Q|^{\frac{\alpha}{nq_{0}^{'}}}}\frac{(w^{p}(F\cap Q))^{\frac{1}{p}}}{(w^{q}(Q))^{\frac{1}{q}}}=[w^{q}]_{A_{\frac{p}{p_{0}},\frac{q}{p_{0}}}^{\mathcal{R}}}^{\frac{1}{p_{0}}}\frac{1}{|Q|^{\frac{\alpha}{nq_{0}^{'}}}}\frac{\left(\int_{Q}fw^{p}dx\right)^{\frac{1}{p}}}{(w^{q}(Q))^{\frac{1}{q}}}.\tag{3.3}
\end{align*}
Besides, apply (1) of Proposition 2.6, we have $\tilde{p}$ that $q<\tilde{p}<q_{0}$ so that $w^{q}\in RH_{\left(q_{0}/\tilde{p}\right)^{'}}$, then
\begin{align*}
\left\langle\chi_{G^{\prime}} w^{q}\right\rangle_{\alpha,q_0^{\prime}, Q}|Q|&=|Q|^{\frac{\alpha}{nq^{'}_{0}}}\langle\chi_{G^{'}}w^{q}\rangle_{q_{0}^{'},Q}|Q|=|Q|^{\frac{\alpha}{nq_{0}^{'}}}|Q|^{1-\frac{1}{q_{0}^{\prime}}}\left(\int_Q\left(\chi_{G^{\prime}} w^{q}\right)^{q_{0}^{'}}dx\right)^{\frac{1}{q_{0}^{'}}}\\
&\leq|Q|^{\frac{\alpha}{nq_{0}^{'}}}|Q|^{1-\frac{1}{q_{0}^{'}}}\left(\int_{Q}w^{q\frac{q_{0}}{q_{0}-\tilde{p}}}dx\right)^{\frac{1}{\tilde{p}}-\frac{1}{q_{0}}}\left(\int_{Q}(\chi_{G^{'}})^{\tilde{p}^{'}}w^{q}dx\right)^{\frac{1}{\tilde{p}^{'}}}\\
&=|Q|^{\frac{\alpha}{nq_{0}^{'}}}\left[\frac{\left(\frac{1}{|Q|}\int_{Q}w^{q\left(\frac{q_{0}}{\tilde{p}}\right)^{'}}dx\right)^{\frac{1}{\left(\frac{q_{0}}{\tilde{p}}\right)^{'}}}}{\frac{1}{|Q|}\int_{Q}w^{q}dx}\right]^{\frac{1}{\tilde{p}}}\left(\frac{1}{w^{q}(Q)}\int_{Q}\chi_{G^{'}}w^{q}dx\right)^{\frac{1}{\tilde{p}^{'}}} w^{q}(Q)\\
&\leq[w^{q}]^{\frac{1}{\tilde{p}}}_{RH_{\left(\frac{q_{0}}{\tilde{p}}\right)^{'}}}|Q|^{\frac{\alpha}{nq_{0}^{'}}}\left(\frac{1}{w^{q}(Q)}\int_{Q}\chi_{G^{'}}w^{q}dx\right)^{\frac{1}{\tilde{p}^{'}}}w^{q}(Q).\tag{3.4}
\end{align*}
Here we have used H\"older's inequality in the second inequality. According to $(3.1)$, $(3.3)$ and $(3.4)$, there exists an $\eta$ sparse collection $\mathcal{S}$, such that
\begin{align*}
    |\langle T f, \chi_{G^{'}}w^{q}\rangle| &\lesssim\sum_{Q \in S}\langle f\rangle_{p_0, Q}\langle \chi_{G^{'}}w^{q}\rangle_{\alpha,q_0^{\prime}, Q}|Q|\\
    &\leq[w^{q}]_{A_{\frac{p}{p_{0}},\frac{q}{p_{0}}}^{\mathcal{R}}}^{\frac{1}{p_{0}}}[w^{q}]^{\frac{1}{\tilde{p}}}_{RH_{\left(\frac{q_{0}}{\tilde{p}}\right)^{'}}}\sum_{Q\in\mathcal{S}}\frac{\left(\int_{Q}fw^{p}dx\right)^{\frac{1}{p}}}{(w^{q}(Q))^{\frac{1}{q}}}\left(\frac{1}{w^{q}(Q)}\int_{Q}\chi_{G^{'}}w^{q}dx\right)^{\frac{1}{\tilde{p}^{'}}} w^{q}(Q)\\
    &=[w^{q}]_{A_{\frac{p}{p_{0}},\frac{q}{p_{0}}}^{\mathcal{R}}}^{\frac{1}{p_{0}}}[w^{q}]^{\frac{1}{\tilde{p}}}_{RH_{\left(\frac{q_{0}}{\tilde{p}}\right)^{'}}}\sum_{Q\in\mathcal{S}}\left(\frac{\int_{Q}fw^{p}dx}{w^{q}(Q)^{1-\frac{p\alpha}{nq_{0}^{'}}}}\right)^{\frac{1}{p}}\left(\frac{1}{w^{q}(Q)}\int_{Q}\chi_{G^{'}}w^{q}dx\right)^{\frac{1}{\tilde{p}^{'}}} w^{q}(Q).\tag{3.5}
\end{align*}
Here we could assume that $G^{'}\cap Q\neq\varnothing$, else $\int_{Q}\chi_{G^{'}}w^{q}dx=0$. For $Q$, $P\in\mathcal{D}$, where $\mathcal{D}$ is the dyadic grid associated to the sparse collection $\mathcal{S}$, we have $Q\cap P\neq\varnothing$, $Q\subseteq P$, or $P\subseteq Q$. As $P\subseteq\Omega$, if $Q\subset P\subseteq\Omega$ for some $P\in\mathcal{P}$, then $Q\cap G^{'}=Q\cap G\cap\Omega^{c}=\varnothing$, which leads to a contradiction. Thus if we have $Q\cap P\neq\varnothing$ for some $P\in\mathcal{P}$, then $P\subseteq Q$, which means $Q\cap\Omega=\bigcup\limits_{\substack{P\subseteq Q\\P\in\mathcal{P}}}P$. Hence, we have
\begin{align*}
    \int_{Q}fw^{p}dx&=\int_{Q}gw^{q}dx+\int_{Q}bw^{q}dx\\
    &=\int_{Q}gw^{q}dx+\int_{Q\cap\Omega}bw^{q}dx\quad(\text{as supp}\ b\subseteq\Omega)\\
    &=\int_{Q}gw^{q}dx+\sum\limits_{\substack{P\subseteq Q\\P\in\mathcal{P}}}\int_{P}bw^{q}dx=\int_{Q}gw^{q}dx.\tag{3.6}
\end{align*}
Thus
\begin{align*}
    \text{RHS of }(3.5)&\leq C_{1}\sum_{Q\in\mathcal{S}}\left( \inf_{E_{Q}}\mathcal{M}^{\mathcal{D}}_{p\alpha/q_{0}^{'},w^{q}}g\right)^{\frac{1}{p}}\left( \inf_{E_{Q}}\mathcal{M}^{\mathcal{D}}_{w^{q}}\chi_{G^{'}}\right)^{\frac{1}{\tilde{p}^{'}}}w^{q}(E_{Q})\\
    &\leq C_{1}\int\left( \mathcal{M}^{\mathcal{D}}_{p\alpha/q_{0}^{'},w^{q}}g\right)^{\frac{1}{p}}\left( \mathcal{M}^{\mathcal{D}}_{w^{q}}\chi_{G^{'}}\right)^{\frac{1}{\tilde{p}^{'}}}w^{q}dx\\
    &\leq C_{1}\left(\int\left(\mathcal{M}^{\mathcal{D}}_{p\alpha/q_{0}^{'},w^{q}}g\right)^{\frac{\theta}{p}}w^{q}dx\right)^{\frac{1}{\theta}}\left(\int\left(\mathcal{M}^{\mathcal{D}}_{w^{q}}\chi_{G^{'}}\right)^{\frac{\theta^{'}}{\tilde{p}^{'}}}w^{q}dx\right)^{\frac{1}{\theta^{'}}}\tag{3.7}
\end{align*}
for any $q<\theta<\tilde{p}$, where $C_{1}=[w^{q}]_{A_{{p}/{p_{0}},{q}/{p_{0}}}^{\mathcal{R}}}^{\frac{1}{p_{0}}}[w^{q}]^{\frac{1}{\tilde{p}}}_{RH_{\left({q_{0}}/{\tilde{p}}\right)^{'}}}D^{\eta}_{w^{q}}.$ Then it follows from Proposition 2.7 that
\begin{align*}
    (3.7)&\leq C_{1}C_{2}||g||^{\frac{1}{p}}_{L^{s}(w^{q})}||\chi_{G^{'}}||^{\frac{1}{\tilde{p}^{'}}}_{L^{\frac{\theta^{'}}{\tilde{p}^{'}}}(w^{q})}\\
    &\leq C_{1}C_{2}\left(||g||_{L^{\infty}}^{\frac{1}{s^{'}}}||g||_{L^{1}(w^{q})}^{\frac{1}{s}}\right)^{\frac{1}{p}}\left(w^{q}(G^{'})\right)^{\frac{1}{\theta^{'}}}\\
    &\lesssim C_{1}C_{2}(D^{2^{n}}_{w^{q}})^{\frac{1}{s^{'}p}}w^{q}(G^{'})^{\frac{1}{\theta^{'}}+\frac{1}{p}\left(\frac{1}{s}-1\right)}=C_{1}C_{2}(D^{2^{n}}_{w^{q}})^{\frac{1}{s^{'}p}}w^{q}(G^{'})^{1-\frac{1}{q}},\tag{3.8}
\end{align*}
where $1/s-\theta/p=p\alpha/(nq_{0}^{'})$ and $C_{2}=\left(1+s^{'}p/\theta\right)^{\left(1-p\alpha/(nq_{0}^{'})\right)/{p}}\left(\theta^{'}/{\tilde{p}^{'}}\right)^{'1/{\tilde{p}^{'}}}$. Denote $r=\left({q_{0}}/{q}\right)^{'}$, we choose $\tilde{p}$ as in Proposition 2.6 (1) that $\left({q_{0}}/{\tilde{p}}\right)^{'}=r+(r-1)/(cr[w^{q}]^{r}_{RH_{r}})$, we denote this as $r+A$ and it holds that $[w^{q}]_{RH_{\left({q_{0}}/{\tilde{p}}\right)^{'}}}\lesssim [w^{q}]_{RH_{r}}$. Let $\theta=q\left(1+{1}/{q}-{1}/{\tilde{p}}\right)$, it is obvious that $\theta>q$, and $\theta=q+(\tilde{p}-q)/{\tilde{p}}<q+\tilde{p}-q=\tilde{p}$. Then
\begin{align*}
    C_{2}&=\left(1+\frac{p}{\theta}\frac{1}{1-\frac{p}{\theta}-\frac{p\alpha}{nq_{0}^{'}}}\right)^{\frac{1}{p}-\frac{\alpha}{nq_{0}^{'}}}\left(\frac{\theta^{'}}{\theta^{'}-\tilde{p}^{'}}\right)^{\frac{1}{\tilde{p}^{'}}}
    =\left(\frac{\theta}{\theta-q}\right)^{\frac{1}{q}}\left(\frac{\theta(\tilde{p}-1)}{\tilde{p}-\theta}\right)^{\frac{1}{\tilde{p}^{'}}}=\left(\frac{q+1-\frac{q}{\tilde{p}}}{1-\frac{q}{\tilde{p}}}\right)^{\frac{1}{q}+\frac{1}{\tilde{p}^{'}}},
\end{align*}
and 
$$\frac{1}{\frac{q_{0}}{q}-\frac{q_{0}}{\tilde{p}}}=\frac{1}{r^{'}-(r+A)^{'}}=r-1+\frac{(r-1)^{2}}{A}=\frac{q}{q_{0}-q}\left(1+cr[w^{q}]_{RH_{r}}^{r}\right).$$ 
Thus
\begin{equation}
    C_{2}=\left(1+\frac{1}{\frac{1}{q}-\frac{1}{\tilde{p}}}\right)^{\frac{1}{q}-\frac{1}{\tilde{p}}+1}\lesssim1+\frac{1}{\frac{1}{q}-\frac{1}{\tilde{p}}}=1+\frac{qq_{0}}{q_{0}-q}\left(1+cr[w^{q}]_{RH_{r}}^{r}\right)\lesssim[w^{q}]_{RH_{r}}^{r}.\tag{3.9}
\end{equation}
Observe that ${1}/{(s^{'}p)}={1}/{q}-{1}/{\theta}$, finally according to $(1)$, $(8)$ and $(9)$, we bound
\begin{align*}
    \|T\|_{L^{p, 1}(w^{p}) \rightarrow L^{q, \infty}(w^{q})}&\lesssim C_{1}C_{2}(D^{2^{n}}_{w^{q}})^{\frac{1}{s^{'}p}}\lesssim [w^{q}]_{A_{\frac{p}{p_{0}},\frac{q}{p_{0}}}^{\mathcal{R}}}^{\frac{1}{p_{0}}}[w^{q}]_{RH_{\left(\frac{q_{0}}{q}\right)^{'}}}^{\left(\frac{q_{0}}{q}\right)^{'}+\frac{1}{q}}D^{\eta}_{w^{q}}\left(D_{w^{q}}^{2^{n}}\right)^{\frac{1}{q}}.
\end{align*}
\underline{The case $\Omega\neq\varnothing\ \text{and}\ q_{0}=\infty$}. With the estimate $(3.3)$,
\begin{align*}
    |\langle T f, \chi_{G^{'}}w^{q}\rangle| &\lesssim\sum_{Q \in S}\langle f\rangle_{p_0, Q}\langle \chi_{G^{'}}w^{q}\rangle_{\alpha,1, Q}|Q|\\
    &\leq [w^{q}]_{A_{\frac{p}{p_{0}},\frac{q}{p_{0}}}^{\mathcal{R}}}^{\frac{1}{p_{0}}}\sum_{Q\in\mathcal{S}}\frac{\left(\int_{Q}fw^{p}dx\right)^{\frac{1}{p}}}{(w^{q}(Q))^{\frac{1}{q}}}\langle \chi_{G^{'}}\rangle_{Q}^{w^{q}}w^{q}(Q)\\
    &\leq [w^{q}]_{A_{\frac{p}{p_{0}},\frac{q}{p_{0}}}^{\mathcal{R}}}^{\frac{1}{p_{0}}}D^{\eta}_{w^{q}}\int\left( \mathcal{M}^{\mathcal{D}}_{p\alpha/q_{0}^{'},w^{q}}g\right)^{\frac{1}{p}} \mathcal{M}^{\mathcal{D}}_{w^{q}}\chi_{G^{'}}w^{q}dx\\   
    &\leq C\left(\int\left(\mathcal{M}^{\mathcal{D}}_{p\alpha/q_{0}^{'},w^{q}}g\right)^{\frac{q+1}{p}}w^{q}dx\right)^{\frac{1}{q+1}}\left(\int\left(\mathcal{M}^{\mathcal{D}}_{w^{q}}\chi_{G^{'}}\right)^{(q+1)^{'}}w^{q}dx\right)^{\frac{1}{(q+1)^{'}}}
\end{align*}
With the same discussion as before, we can prove
$$ \|T\|_{L^{p, 1}(w^{p}) \rightarrow L^{q, \infty}(w^{q})}\lesssim[w^{q}]_{A_{\frac{p}{p_{0}},\frac{q}{p_{0}}}^{\mathcal{R}}}^{\frac{1}{p_{0}}}D^{\eta}_{w^{q}}\left(D_{w^{q}}^{2^{n}}\right)^{\frac{1}{q}-\frac{1}{q+1}}.$$
\underline{The case $\Omega=\varnothing$}. The fact $(3.2)$ asserts that for almost every $x\in\mathbf{R}^{n}$, $fw^{p-q}\lesssim w^{q}(G)^{-1}$, thus there is no need to involve the Calder\'{o}n-Zygmund decomposition. Instead, we replace the good part $g$ in the proof above with function $fw^{p-q}$, and this leads to the same bound.\\
\indent Finally, we estimate the constant $D^{\eta}_{w^{q}}$ as following: according to H\"older's inequality, for a cube $Q$ and subset $E\subseteq Q$ that $|E|\geq\eta|Q|$, 
$$\left(\frac{|E|}{|Q|^{1-\frac{p_{0}\alpha}{nq_{0}^{'}}}}\right)^{\frac{1}{p_{0}}}\leq[w^{q}]_{A_{\frac{p}{p_{0}},\frac{q}{p_{0}}}^{\mathcal{R}}}^{\frac{1}{p_{0}}}\frac{\left(w^{p}(E)\right)^{\frac{1}{p}}}{\left(w^{q}(Q)\right)^{\frac{1}{q}}}\leq[w^{q}]_{A_{\frac{p}{p_{0}},\frac{q}{p_{0}}}^{\mathcal{R}}}^{\frac{1}{p_{0}}}\frac{\left(w^{q}(E)\right)^{\frac{1}{q}}}{\left(w^{q}(Q)\right)^{\frac{1}{q}}}|E|^{\frac{1}{p}-\frac{1}{q}},$$
where ${1}/{p}-{1}/{q}={\alpha}/{nq_{0}^{'}}$. Hence
$$\frac{w^{q}(Q)}{w^{q}(E)}\leq[w^{q}]_{A_{\frac{p}{p_{0}},\frac{q}{p_{0}}}^{\mathcal{R}}}^{\frac{q}{p_{0}}}\eta^{-\frac{\alpha q}{nq_{0}^{'}}+\frac{q}{p_{0}}}.\ \square$$
\subsection{Proof of Theorem 1.3}
\indent With a similar approach as before, we assume $||f||_{L^{p}}=1$ and use the equivalence 
\begin{equation}
    ||wT(fw^{-1})||_{L^{q,\infty}}\sim\sup\limits_{0<|G|<\infty}\inf_{\substack{G^{'}\subseteq G\\|G^{'}|\geq\frac{1}{2}|G|}}|G^{'}|^{-1+\frac{1}{q}}|\langle T(fw^{-1}), w\chi_{G^{'}}\rangle|.\tag{3.10}
\end{equation}
For fixed $G$ with $0<|G|<\infty$, let
$$\Omega=\left\{\mathcal{M}^{\mathcal{D}}(|f|^{p})>\frac{2}{|G|}\right\},$$
and let $\mathcal{P}$ be the family of those maximal cubes in $\mathcal{D}$ that $\langle|f|^{p}\rangle_{Q}>{2}/{|G|}$.\\
\underline{The case $\Omega\neq\varnothing$}. Then $\bigcup_{P\in\mathcal{P}}P=\Omega$ and according to Proposition 2.7 it holds that $|\Omega|\leq|G|\int|f|^{p}dx/2={|G|}/{2}$. Let $G^{'}=G\cap\Omega^{c}$, we have $|G^{'}|\geq|G|/2$. Next consider the Calder\'{o}n-Zygmund decomposition that $|f|^{p}=g+b$, where $g=\sum_{P\in\mathcal{P}}\langle|f|^{p}\rangle_{P}\chi_{P}+|f|^{p}\chi_{\Omega^{c}}$, $b=\sum_{P\in\mathcal{P}}b_{P}$ that $b_{P}=\left(|f|^{p}-\langle|f|^{p}\rangle_{P}\right)\chi_{P}$. We know that $||g||_{L^{1}}=||f||_{L^{p}}^{p}=1$, $||g||_{L^{\infty}}\lesssim|G|^{-1}$, supp$(b)\subseteq\Omega$ and the bad part $b_{P}$ enjoys the cancelling property that for any $P\in\mathcal{P}$, $\langle b_{P}\rangle_{P}=0$.\\
\indent Denote $s=\left({q_{0}}/{q}\right)^{'}$. Since $w^{q}\in RH_{s}$, Proposition 2.6 asserts there exists a constant $v>s$, such that $w^{q}\in RH_{v}$. When $q_{0}=\infty$ i.e. $s=1$, $v$ could be taken as $v=1+{d}([w^{q}]_{A_{\infty}})^{-1}$ that $[w^{q}]_{RH_{v}}\lesssim1$. While when $q_{0}<\infty$, we take $v$ as $v=s+(s-1)/(cs[w^{q}]_{RH_{s}}^{s})$, there it holds that $[w^{q}]_{RH_{v}}\leq[w^{q}]_{RH_{s}}$. Next we define $\ell>0$ by ${1}/{\ell}={1}/{q_{0}^{'}}-{1}/{(qv)}$. Here 
$$\frac{1}{q_{0}^{'}}>\frac{1}{q_{0}^{'}}-\frac{1}{qv}>\frac{1}{q_{0}^{'}}-\frac{1}{qs}=\frac{1}{q_{0}^{'}}-\frac{q_{0}-q}{qq_{0}}=1-\frac{1}{q}=\frac{1}{q^{'}},$$
thus $q_{0}^{'}<\ell<q^{'}$. Define $r={1}/{q}+{1}/{\ell}>{1}/{q}+{1}/{q^{'}}=1$, we have ${(qr)^{'}}/{c}=r$. We first consider the case $p>p_{0}$:
\begin{align*}
    |\langle T(fw^{-1}), w\chi_{G^{'}}\rangle|&\leq\sum_{Q\in\mathcal{S}}\langle fw^{-1}\rangle_{p_{0},Q}\langle w\chi_{G^{'}}\rangle_{q_{0}^{'},Q}|Q|^{1+\frac{\alpha}{nq_{0}^{'}}}\\
    &\leq\sum_{Q\in\mathcal{S}}\langle f\rangle_{p,Q}\langle w^{-1}\rangle_{\frac{1}{\frac{1}{p_{0}}-\frac{1}{p}},Q}\langle w\rangle_{qv,Q}\langle \chi_{G^{'}}\rangle_{\ell,Q}|Q|^{1+\frac{\alpha}{nq_{0}^{'}}}\\
    &\leq[w^{q}]_{A_{\left(\frac{1}{p_{0}}-\frac{1}{p}\right)q+1}}^{\frac{1}{q}}\sum_{Q\in\mathcal{S}}\langle f\rangle_{p,Q}\langle w\rangle_{q,Q}^{-1}\langle w\rangle_{qv,Q}\langle \chi_{G^{'}}\rangle_{\ell,Q}|Q|^{1+\frac{\alpha}{nq_{0}^{'}}}\tag{3.11}\\
    &\leq[w^{q}]_{A_{\left(\frac{1}{p_{0}}-\frac{1}{p}\right)q+1}}^{\frac{1}{q}}[w^{q}]_{RH_{v}}^{\frac{1}{q}}\sum_{Q\in\mathcal{S}}\langle f\rangle_{p,Q}\langle\chi_{G^{'}}\rangle_{\ell,Q}|Q|^{1+\frac{\alpha}{nq_{0}^{'}}}.\tag{3.12}
\end{align*}
Here in the first inequality, we have used H\"older's inequality, and the second and the third one due to the definition of corresponding weight. For the case $p=p_{0}$, we instead use the estimate $\langle fw^{-1}\rangle_{p_{0},Q}\leq\langle f\rangle_{p_{0},Q}||w^{-q}||_{L^{\infty}(Q)}^{\frac{1}{q}}\leq[w^{q}]_{A_{1}}^{\frac{1}{q}}\langle f\rangle_{p_{0},Q}\langle w\rangle_{q,Q}^{-1}$ in the second inequality. \\
Then with a discussion similar to (3.6), we can check the cancelling property of $b$ as
$$\int_{Q}|f|^{p}dx=\int_{Q}gdx,\ \ \text{for any }Q\in\mathcal{S}\ \text{that }Q\cap G^{'}\neq\varnothing.$$
Thus according to Proposition 2.7 we have
\begin{align*}
    (3.12)&\lesssim[w^{q}]_{A_{\left(\frac{1}{p_{0}}-\frac{1}{p}\right)q+1}}^{\frac{1}{q}}[w^{q}]_{RH_{v}}^{\frac{1}{q}}\sum_{Q\in\mathcal{S}}\left(\inf_{E_{Q}}\mathcal{M}^{\mathcal{D}}_{\frac{p\alpha}{q_{0}^{'}}}g\right)^{\frac{1}{p}}\left(\inf_{E_{Q}}\mathcal{M}^{\mathcal{D}}\chi_{G^{'}}\right)^{\frac{1}{\ell}}|E_{Q}|\\
    &\leq[w^{q}]_{A_{\left(\frac{1}{p_{0}}-\frac{1}{p}\right)q+1}}^{\frac{1}{q}}[w^{q}]_{RH_{v}}^{\frac{1}{q}}\int\left(\mathcal{M}^{\mathcal{D}}_{\frac{p\alpha}{q_{0}^{'}}}g\right)^{\frac{1}{p}}\left(\mathcal{M}^{\mathcal{D}}\chi_{G^{'}}\right)^{\frac{1}{\ell}}dx\\
    &\leq[w^{q}]_{A_{\left(\frac{1}{p_{0}}-\frac{1}{p}\right)q+1}}^{\frac{1}{q}}[w^{q}]_{RH_{v}}^{\frac{1}{q}}||\mathcal{M}^{\mathcal{D}}_{\frac{p\alpha}{q_{0}^{'}}}g||_{L^{\frac{qr}{p}}}^{\frac{1}{p}}||\mathcal{M}^{\mathcal{D}}\chi_{G^{'}}||_{L^{r}}^{\frac{1}{\ell}}\\
    &\leq[w^{q}]_{A_{\left(\frac{1}{p_{0}}-\frac{1}{p}\right)q+1}}^{\frac{1}{q}}[w^{q}]_{RH_{v}}^{\frac{1}{q}}\left[\left(1+\frac{t^{'}p}{qr}\right)^{1-\frac{p\alpha}{nq_{0}^{'}}}\right]^{\frac{1}{p}}(r^{'})^{\frac{1}{\ell}}||g||_{L^{t}}^{\frac{1}{p}}||\chi_{G^{'}}||_{L^{r}}^{\frac{1}{\ell}},\tag{3.13}
\end{align*}
where ${1}/{t}-{p}/{(qr)}={p\alpha}/{(nq_{0}^{'})}$, ${p}/{(qr)}+{p\alpha}/{(nq_{0}^{'})}={p}/{(qr)}+1-{p}/{q}<1$, thus $t>1$. Here 
\begin{equation}
    ||g||_{L^{t}}^{\frac{1}{p}}||\chi_{G^{'}}||_{L^{r}}^{\frac{1}{\ell}}\leq\left(||g||_{L^{1}}^{\frac{1}{t}}||g||_{L^{\infty}}^{\frac{1}{t^{'}}}\right)^{\frac{1}{p}}|G|^{\frac{1}{\ell r}}\leq|G^{'}|^{\frac{1}{\ell r}-\frac{1}{pt^{'}}}=|G^{'}|^{\frac{1}{cr}-\frac{1}{p}+\frac{1}{qr}+\frac{1}{p}-\frac{1}{q}}=|G^{'}|^{1-\frac{1}{q}}.\tag{3.14}
\end{equation}
We then estimate the constant: $1+{t^{'}p}/{(qr)}=t^{'}\left({1}/{t^{'}}+{1}/{t}-{p\alpha}/{(nq_{0}^{'})}\right)=t^{'}\left(1-{p\alpha}/{(nq_{0}^{'})}\right)$, $1-{p\alpha}/{(nq_{0}^{'})}=1-p\left({1}/{p}-{1}/{q}\right)={p}/{q}$, ${1}/{t^{'}}=1-{p}/{(qr)}-{p\alpha}/{(nq_{0}^{'})}={p}/{q}-{p}/{(qr)}={p}/{(qr^{'})}$. Thus
$$\left[\left(1+\frac{t^{'}p}{qr}\right)^{1-\frac{p\alpha}{nq_{0}^{'}}}\right]^{\frac{1}{p}}(r^{'})^{\frac{1}{\ell}}=(r^{'})^{\frac{1}{q}+\frac{1}{\ell}}=(r^{'})^{r},$$
and $r={1}/{q_{0}^{'}}+{1}/{(qv^{'})}$. For the case $s=1$ i.e. $q_{0}=\infty$, it holds that $v^{'}\sim[w^{q}]_{A_{\infty}}$, $r=1+{1}/{(qv{'})}$, $r^{'}={r}/{(r-1)}=1+qv^{'}$, hence
\begin{equation}
    (r^{'})^{r}=(1+qv^{'})^{1+\frac{1}{qv^{'}}}\lesssim1+qv^{'}\lesssim[w^{q}]_{A_{\infty}}.\tag{3.15}
\end{equation}
Together with $(3.10),(3.13)$ it follows that $||wT(w^{-1}\cdot)||_{L^{p}\rightarrow L^{q,\infty}}\lesssim[w^{q}]_{A_{\left(\frac{1}{p_{0}}-\frac{1}{p}\right)q+1}}^{\frac{1}{q}}[w^{q}]_{A_{\infty}}$.\\
Else if $s>1$ i.e. $q_{0}<\infty$, we have
\begin{align*}
    (r^{'})^{r}=\left(\frac{\frac{1}{q_{0}^{'}}+\frac{1}{qv^{'}}}{\frac{1}{qv^{'}}-\frac{1}{q_{0}}}\right)^{\frac{1}{q_{0}^{'}}+\frac{1}{qv^{'}}}\lesssim\left(\frac{1}{\frac{1}{qv^{'}}-\frac{1}{q_{0}}}\right)^{\frac{1}{q_{0}^{'}}+\frac{1}{qv^{'}}}\lesssim\left(\frac{1}{1-\frac{v^{'}}{s^{'}}}\right)^{\frac{1}{q_{0}^{'}}+\frac{1}{qv^{'}}}.
\end{align*}
As $v=s+{(s-1)}/{(cs[w^{q}]_{RH_{s}}^{s})}$, we have ${(1-{v^{'}}/{s^{'}})}^{-1}={s(v-1)}/{(v-s)}\lesssim[w^{q}]_{RH_{s}}^{s}$, thus $(r^{'})^{r}\lesssim[w^{q}]_{RH_{s}}^{s\left(\frac{1}{q_{0}^{'}}+\frac{1}{q}\right)}=[w^{q}]_{RH_{\left({q_{0}}/{q}\right)^{'}}}^{\left(\frac{q_{0}}{q}\right)^{'}+\frac{1}{q}}$. Together with $(10),(13)$ and $[w^{q}]_{RH_{v}}\lesssim[w^{q}]_{RH_{s}}$, it follows that 
$$||wT(w^{-1}\cdot)||_{L^{p}\rightarrow L^{q},\infty}\lesssim[w^{q}]_{A_{\left(\frac{1}{p_{0}}-\frac{1}{p}\right)q+1}}^{\frac{1}{q}}[w^{q}]_{RH_{\left(\frac{q_{0}}{q}\right)^{'}}}^{\left(\frac{q_{0}}{q}\right)^{'}+\frac{2}{q}}.$$
\underline{The case $\Omega=\varnothing$}. With a similar discussion as in the proof of Theorem 1.1, we replace $g$ with $f$ in the proof above, then almost everywhere $f\lesssim|G|^{-1}$, which eventually leads to the same bound. $\square$\\
~\\
\textbf{Acknowledgement}\\
\indent The author would like to extend sincere gratitude to his supervisor, Professor Y. Tsutsui, who instructs the
 author in the weight theory and gives lots of valuable suggestions.







\end{document}